\documentclass{article}%
\usepackage{amsfonts}
\usepackage{amsmath}%
\usepackage{amssymb}%
\usepackage{graphicx}
\providecommand{\U}[1]{\protect\rule{.1in}{.1in}}

\begin{document}

\title{Evaluating Polynomials Over the Unit Disk \\and the Unit Ball}
\author{Kendall Atkinson, University of Iowa, Iowa City, Iowa
\and Olaf Hansen and David Chien, California State University, San Marcos, CA}
\maketitle

\begin{abstract}
We investigate the use of orthonormal polynomials over the unit disk
$\mathbb{B}_{2}$ in $\mathbb{R}^{2}$ and the unit ball $\mathbb{B}_{3}$ in
$\mathbb{R}^{3}$. An efficient evaluation of an orthonormal polynomial basis
is given, and it is used in evaluating general polynomials over $\mathbb{B}%
_{2}$ and $\mathbb{B}_{3}$. The least squares approximation of a function $f$
on the unit disk by polynomials of a given degree is investigated, including
how to write a polynomial using the orthonormal basis. \textsc{Matlab} codes
are given.

\end{abstract}

\section{Introduction\label{intro}}

A standard way to write a multivariate polynomial of degree $n$ over
$\mathbb{R}^{2}$ is%
\[
p\left(  x,y\right)  =\sum_{j=0}^{n}\sum_{k=0}^{j}a_{j,k}x^{j}y^{j-k}.
\]
The space of all such polynomials is denoted by $\Pi_{n}$. We consider here
the alternative formulation%
\begin{equation}
p\left(  x,y\right)  =\sum_{j=0}^{n}\sum_{k=0}^{j}b_{j,k}\varphi_{j,k}\left(
x,y\right)  \label{eq1}%
\end{equation}
with $\left\{  \varphi_{j,k}\mid0\leq k\leq j,\ 0\leq j\leq n\right\}  $ an
orthonormal basis of the set of $\Pi_{n}$ over the closed unit disk
$\mathbb{B}_{2}$, for each $n\geq0$. There is a large literature on such
orthonormal polynomials; and in contrast to the univariate case, there are
many possible choices for this basis. See Dunkl and Xu \cite{DX} and Xu
\cite{Xu2004} for an investigation of such multivariate orthonormal
polynomials and a number of particular examples.

To use (\ref{eq1}), it is important to be able to evaluate the orthonormal
polynomials $\left\{  \varphi_{j,k}\right\}  $ efficiently, just as is true
with univariate polynomials. We consider a particularly good set of such
polynomials in Section \ref{sec2}, one that seems much superior to other
choices. In the univariate case, the best choices are based on using the
triple recursion relation of the particular family $\left\{  \varphi
_{n}\right\}  $ being used. This extends to the multivariate case. We
investigate a particular choice of \ an orthonormal basis for $\Pi_{n}$ that
leads to an efficient way to evaluate the expression (\ref{eq1}) by making use
of the triple recursion relation it satisfies. Following that, in Section
\ref{sec3}, we also consider the calculation of the least squares
approximation over $\Pi_{n}$ of a given function $f\left(  x,y\right)  $. In
Section \ref{sec4}, these results are extended to polynomials over the unit
ball. Finally, in Section \ref{sec5}, \textsc{Matlab} codes are given for all
of the problems being discussed.

\section{Evaluating an orthonormal polynomial basis\label{sec2}}

We review some notation and results from Dunkl and Xu \cite{DX} and Xu
\cite{Xu2004}. For convenience, we initially denote a point in the unit disk
by $x=\left(  x_{1},x_{2}\right)  $, and later we revert to the more standard
use of $\left(  x,y\right)  $. We consider only the standard $L^{2}$ inner
product
\begin{equation}
\left(  p,q\right)  =\int_{\mathbb{B}_{2}}p\left(  x\right)  q\left(
x\right)  \,dx. \label{inner}%
\end{equation}
Define
\[
\mathcal{V}_{n}=\left\{  p\in\Pi_{n}\mid\left(  p,q\right)  =0,\ \forall
q\in\Pi_{n-1}\right\}  ,\quad n\geq1,
\]
and let $\mathcal{V}_{0}$ denote the one dimensional space of constant
functions. Thus
\[
\Pi_{n}=\mathcal{V}_{0}\oplus\cdots\oplus\mathcal{V}_{n}%
\]
is an orthogonal decomposition of $\Pi_{n}$. It is standard to give an
orthonormal basis for each space $\mathcal{V}_{n}$ as the way to give an
orthonormal basis of $\Pi_{n}$. The dimension of $\mathcal{V}_{n}$ equals
$n+1$, and the dimension of $\Pi_{n}$ equals%
\begin{equation}
N_{n}=\frac{1}{2}\left(  n+1\right)  \left(  n+2\right)  . \label{dimen}%
\end{equation}

Introduce
\[
\mathbb{P}_{n}\mathbb{=}\left[  Q_{n}^{0},Q_{n}^{1},\dots,Q_{n}^{n}\right]
^{\text{T}},\quad\quad n\geq0,
\]
with $\left\{  Q_{n}^{0},Q_{n}^{1},\dots,Q_{n}^{n}\right\}  $ an orthonormal
basis of $\mathcal{V}_{m}$. The triple recursion relation for $\left\{
\mathbb{P}_{m}\right\}  $ is given by%
\begin{equation}
x_{i}\mathbb{P}_{n}\left(  x\right)  =A_{n,i}\mathbb{P}_{n+1}\left(  x\right)
+B_{n,i}\mathbb{P}_{n}\left(  x\right)  +A_{n-1,i}^{\text{T}}\mathbb{P}%
_{n-1}\left(  x\right)  ,\quad\quad i=1,2,\quad n\geq1 \label{triple}%
\end{equation}
The matrices $A_{n,i}$ and $B_{n,i}$ are $\left(  n+1\right)  \times\left(
n+2\right)  $ and $\left(  n+1\right)  \times\left(  n+1\right)  $,
respectively, and they are defined as follows:%
\begin{align*}
A_{n,i}  &  =\int_{\mathbb{B}_{2}}x_{i}\mathbb{P}_{n}\left(  x\right)
\mathbb{P}_{n+1}^{\text{T}}\left(  x\right)  \,dx\medskip\\
B_{n,i}  &  =\int_{\mathbb{B}_{2}}x_{i}\mathbb{P}_{n}\left(  x\right)
\mathbb{P}_{n}^{\text{T}}\left(  x\right)  \,dx
\end{align*}
For additional details, see Xu \cite[Thm. 2.1]{Xu2004}. One wants to use the
relation (\ref{triple}) to solve for $\mathbb{P}_{n+1}\left(  x\right)  $.
This amounts to solving an overdetermined system of $2\left(  n+1\right)  $
equations for the $n+2$ components of $\mathbb{P}_{n+1}\left(  x\right)  $.
The expense of this will depend on the structure of the matrices $A_{n,i}$ and
$B_{n,i}$. There is a well-known choice that leads, fortunately, to the
matrices $B_{n,i}$ being zero and the matrices $A_{n,i}$ being very sparse.

To define this choice, begin by recalling the Gegenbauer polynomials $\left\{
C_{n}^{\lambda}\left(  t\right)  \right\}  $. They can be obtained using the
following generating function:%
\[
\left(  1-2rt+r^{2}\right)  ^{-\lambda}=\sum_{n=0}^{\infty}C_{n}^{\lambda
}\left(  t\right)  r^{n},\quad\quad\left\vert r\right\vert <1,\quad\left\vert
t\right\vert \leq1
\]
For particular cases,%
\[
C_{0}^{\lambda}\left(  t\right)  \equiv1,\quad\quad C_{1}^{\lambda}\left(
t\right)  =2\lambda t,\quad\quad C_{2}^{\lambda}\left(  t\right)
=\lambda\left(  2\left(  \lambda+1\right)  t^{2}-1\right)  ,
\]%
\[
C_{3}^{\lambda}\left(  t\right)  =\frac{2}{3}\lambda\left(  \lambda+1\right)
t\left(  \left(  2\lambda+4\right)  t^{2}-3\right)  .
\]
Their triple recursion relation is given by%
\[
C_{n+1}^{\lambda}\left(  t\right)  =\frac{2\left(  n+\lambda\right)  }%
{n+1}tC_{n}^{\lambda}\left(  t\right)  -\frac{n+2\lambda-1}{n+1}%
C_{n-1}^{\lambda}\left(  t\right)  ,\quad\quad n\geq1.
\]
These polynomials are orthogonal over $\left(  -1,1\right)  $ with respect to
the inner product%
\[
\left(  f,g\right)  =\int_{-1}^{1}\left(  1-t^{2}\right)  ^{\lambda-\frac
{1}{2}}f\left(  t\right)  g\left(  t\right)  \,dt,
\]
and for $\lambda=\frac{1}{2}$ they are the Legendre polynomials. For
additional information on the Gegenbauer polynomials, see \cite[Chap.
18]{nist}.

Return to the use of $\left(  x,y\right)  $ in place of $\left(  x_{1}%
,x_{2}\right)  $. Using the Gegenbauer polynomials, introduce%

\begin{equation}
Q_{n}^{k}\left(  x,y\right)  =\frac{1}{h_{k,n}}C_{n-k}^{k+1}\left(  x\right)
\left(  1-x^{2}\right)  ^{\frac{k}{2}}C_{k}^{\frac{1}{2}}\left(  \frac
{y}{\sqrt{1-x^{2}}}\right)  ,\quad\quad\left(  x,y\right)  \in\mathbb{B}_{2},
\label{OrthoBasis}%
\end{equation}
for $n=0,1.\dots$ and $k=0,1,\dots,n$. See Dunkl and Xu \cite[p. 88]{DX}. Note
that%
\[
x^{2}+y^{2}<1\quad\Longrightarrow\quad\frac{\left\vert y\right\vert }%
{\sqrt{1-x^{2}}}<1
\]
The lead constant $h_{k,n}$ is given by%
\[
h_{k,n}^{2}=\frac{\pi}{4^{k}}\frac{\left(  n+k+1\right)  !}{\left(
n+1\right)  \left(  2k+1\right)  \left(  k!\right)  ^{2}\left(  n-k\right)  !}%
\]
and $h_{0,0}^{2}=\pi$. The set $\left\{  Q_{m}^{k}\mid0\leq k\leq m\right\}  $
is an orthonormal basis of $\mathcal{V}_{m}$, and $\left\{  Q_{m}^{k}\mid0\leq
k\leq m,\ 0\leq m\leq n\right\}  $ is an orthonormal basis of $\Pi_{n}$, using
the inner product of (\ref{inner}). Here are the $Q_{m}^{k}$ of degrees
0,1,2,3.%
\begin{equation}
Q_{0}^{0}\left(  x,y\right)  =\frac{1}{\sqrt{\pi}},\quad Q_{1}^{0}\left(
x,y\right)  =\frac{2x}{\sqrt{\pi}},\quad Q_{1}^{1}\left(  x,y\right)
=\frac{2y}{\sqrt{\pi}} \label{poly1}%
\end{equation}%
\begin{equation}
Q_{2}^{0}\left(  x,y\right)  =\frac{1}{\sqrt{\pi}}\left(  4x^{2}-1\right)
,\quad Q_{2}^{1}\left(  x,y\right)  =\sqrt{\frac{24}{\pi}}xy,\quad Q_{2}%
^{2}\left(  x,y\right)  =\sqrt{\frac{2}{\pi}}\left(  3y^{2}+x^{2}-1\right)
\label{poly2}%
\end{equation}%
\begin{equation}%
\begin{tabular}
[c]{ll}%
$Q_{3}^{0}\left(  x,y\right)  =\dfrac{4}{\sqrt{\pi}}x\left(  2x^{2}-1\right)
$ & $\quad Q_{3}^{1}\left(  x,y\right)  =\dfrac{4}{\sqrt{5\pi}}y\left(
6x^{2}-1\right)  ,\smallskip$\\
$Q_{3}^{2}\left(  x,y\right)  =\dfrac{4}{\sqrt{\pi}}x\left(  3y^{2}%
+x^{2}-1\right)  $ & $\quad Q_{3}^{3}\left(  x,y\right)  =\dfrac{4}{\sqrt
{5\pi}}y\left(  5y^{2}-3+3x^{2}\right)  .$%
\end{tabular}
\label{poly3}%
\end{equation}
Because the formula (\ref{OrthoBasis}) is not well-defined at $x=\pm1$, we use%
\[
\lim_{\left(  x,y\right)  \rightarrow\left(  \pm1,0\right)  }\left(
1-x^{2}\right)  ^{\frac{k}{2}}C_{k}^{\frac{1}{2}}\left(  \frac{y}%
{\sqrt{1-x^{2}}}\right)  =\left\{
\begin{array}
[c]{c}%
0,\quad\quad k>0\\
1,\quad\quad k=0
\end{array}
\right.
\]
when evaluating (\ref{OrthoBasis}).

Applying (\ref{triple}) to this choice of orthonormal polynomials leads to
\begin{equation}
x_{i}\mathbb{P}_{n}\left(  x_{1},x_{2}\right)  =A_{n,i}\mathbb{P}_{n+1}\left(
x_{1},x_{2}\right)  +A_{n-1,i}^{\text{T}}\mathbb{P}_{n-1}\left(  x_{1}%
,x_{2}\right)  ,\quad\quad i=1,2,\quad n\geq1. \label{triple2}%
\end{equation}
The coefficient matrices are given by%
\[
A_{n,1}=\left[
\begin{array}
[c]{ccccc}%
a_{0,n} & 0 & \cdots & 0 & 0\\
0 & a_{1,n} &  & 0 & 0\\
\vdots &  & \ddots & \vdots & \vdots\\
0 & 0 & \cdots & a_{n,n} & 0
\end{array}
\right]
\]%
\[
A_{n,2}=\left[
\begin{array}
[c]{cccccc}%
0 & d_{0,n} & 0 & \cdots & 0 & 0\\
c_{1,n} & 0 & d_{1,n} & \ddots & 0 & 0\\
\vdots & \ddots & \ddots & \ddots & \vdots & \vdots\\
0 & \cdots & c_{n-1,n} & 0 & d_{n-1,n} & 0\\
0 & \cdots & 0 & c_{n,n} & 0 & d_{n,n}%
\end{array}
\right]
\]%
\[%
\begin{tabular}
[c]{l}%
$a_{k,n}=\dfrac{1}{2}\sqrt{\dfrac{\left(  n-k+1\right)  \left(  n+k+2\right)
}{\left(  n+1\right)  \left(  n+2\right)  }},\medskip$\\
$d_{k,n}=\dfrac{k+1}{2}\sqrt{\dfrac{\left(  n+k+3\right)  \left(
n+k+2\right)  }{\left(  2k+1\right)  \left(  2k+3\right)  \left(  n+1\right)
\left(  n+2\right)  }},\medskip$\\
$c_{k,n}=-\dfrac{k}{2}\sqrt{\dfrac{\left(  n-k+1\right)  \left(  n-k+2\right)
}{\left(  n+1\right)  \left(  n+2\right)  \left(  2k-1\right)  \left(
2k+1\right)  }}. $%
\end{tabular}
\]
These results are taken from Dunkl and Xu \cite[p. 88]{DX} (in the formula for
$c_{k,n}$, change $n+k+1$ to $n-k+1$).

From the first triple recursion relation in (\ref{triple2}),%

\begin{align*}
x_{1}\left[
\begin{array}
[c]{c}%
Q_{n}^{0}\\
Q_{n}^{1}\\
\vdots\\
Q_{n}^{n}%
\end{array}
\right]   &  =\left[
\begin{array}
[c]{ccccc}%
a_{0,n} & 0 & \cdots & 0 & 0\\
0 & a_{1,n} &  & 0 & 0\\
\vdots &  & \ddots & \vdots & \vdots\\
0 & 0 & \cdots & a_{n,n} & 0
\end{array}
\right]  \left[
\begin{array}
[c]{c}%
Q_{n+1}^{0}\\
Q_{n+1}^{1}\\
\vdots\\
Q_{n+1}^{n}\\
Q_{n+1}^{n+1}%
\end{array}
\right]  \smallskip\\
&  +\left[
\begin{array}
[c]{cccc}%
a_{0,n-1} & 0 & \cdots & 0\\
0 & a_{1,n-1} &  & 0\\
\vdots &  & \ddots & \vdots\\
0 &  &  & a_{n-1,n-1}\\
0 & 0 & \cdots & 0
\end{array}
\right]  \left[
\begin{array}
[c]{c}%
Q_{n-1}^{0}\\
Q_{n-1}^{1}\\
\vdots\\
Q_{n-1}^{n-1}%
\end{array}
\right]
\end{align*}%
\begin{align*}
x_{1}Q_{n}^{i}  &  =a_{i,n}Q_{n+1}^{i}+a_{i,n-1}Q_{n-1}^{i},\quad\quad
i=0,1,\dots,n-1\smallskip\\
x_{1}Q_{n}^{n}  &  =a_{n,n}Q_{n+1}^{n}%
\end{align*}
This allows us to solve for $\left\{  Q_{n+1}^{0},\dots,Q_{n+1}^{n}\right\}
$. The second triple recursion relation in (\ref{triple2}) yields%
\begin{align*}
x_{2}\left[
\begin{array}
[c]{c}%
Q_{n}^{0}\\
Q_{n}^{1}\\
\vdots\\
Q_{n}^{n}%
\end{array}
\right]   &  =\left[
\begin{array}
[c]{cccccc}%
0 & d_{0,n} & 0 & \cdots & 0 & 0\\
c_{1,n} & 0 & d_{1,n} & \ddots & 0 & 0\\
\vdots & \ddots & \ddots & \ddots & \vdots & \vdots\\
0 & \cdots & c_{n-1,n} & 0 & d_{n-1,n} & 0\\
0 & \cdots & 0 & c_{n,n} & 0 & d_{n,n}%
\end{array}
\right]  \left[
\begin{array}
[c]{c}%
Q_{n+1}^{0}\\
Q_{n+1}^{1}\\
\vdots\\
Q_{n+1}^{n}\\
Q_{n+1}^{n+1}%
\end{array}
\right]  \smallskip\\
&  +\left[
\begin{array}
[c]{cccccc}%
0 & c_{1,n-1} & 0 & \cdots &  & 0\\
d_{0,n-1} & 0 & c_{2,n-1} & 0 & \cdots & 0\\
0 & d_{1,n-1} & 0 & c_{3,n-1} &  & 0\\
0 & 0 & d_{2,n-1} & \ddots & \ddots & \\
\vdots &  &  & \ddots & \ddots & c_{n-1,n-1}\\
0 &  &  &  & d_{n-2,n-1} & 0\\
0 & 0 & 0 & \cdots & 0 & d_{n-1,n-1}%
\end{array}
\right]  \left[
\begin{array}
[c]{c}%
Q_{n-1}^{0}\\
Q_{n-1}^{1}\\
\vdots\\
Q_{n-1}^{n-1}%
\end{array}
\right]
\end{align*}
Its last equation is%
\[
x_{2}Q_{n}^{n}=c_{n,n}Q_{n+1}^{n-1}+d_{n,n}Q_{n+1}^{n+1}+d_{n-1,n-1}%
Q_{n-1}^{n-1}%
\]
and from it we can calculate $Q_{n+1}^{n+1}$. Thus,%
\begin{align}
Q_{n+1}^{i}  &  =\frac{x_{1}Q_{n}^{i}-a_{i,n-1}Q_{n-1}^{i}}{a_{i,n}}%
,\quad\quad i=0,1,\dots,n-1\smallskip\label{q1}\\
Q_{n+1}^{n}  &  =\frac{x_{1}Q_{n}^{n}}{a_{n,n}}\smallskip\label{q2}\\
Q_{n+1}^{n+1}  &  =\frac{x_{2}Q_{n}^{n}-c_{n,n}Q_{n+1}^{n-1}-d_{n-1,n-1}%
Q_{n-1}^{n-1}}{d_{n,n}} \label{q3}%
\end{align}

\subsection{Computational cost}

What is the cost of using this to evaluate the orthonormal basis
\[
\mathcal{B}_{n}\equiv\left\{  Q_{m}^{k}\mid0\leq k\leq m,\ 0\leq m\leq
n\right\}  ?
\]
Assume the coefficients $\left\{  a_{i,n},c_{i,n},d_{i,n}\right\}  $ have been
computed. Apply (\ref{q1})-(\ref{q3}) to the computation of $\left\{
Q_{m}^{k}\mid0\leq k\leq m\right\}  $, assuming the lower degree polynomials
of degrees $m-1$ and $m-2$ are known. This requires $4\left(  m+1\right)  $
arithmetic operations. The evaluation of $\{Q_{0}^{0},Q_{1}^{0},Q_{1}^{1}\}$
from (\ref{poly1}) requires 2 arithmetic operations for each choice of
$\left(  x,y\right)  =\left(  x_{1},x_{2}\right)  $. Thus the calculation of
$\mathcal{B}_{n}$ requires%
\begin{equation}
2+4(3+4+\cdots+(n+1))=2\left(  n^{2}+3n-3\right)  \label{BasisCost}%
\end{equation}
arithmetic operations. Recall (\ref{dimen}) that the dimension of $\Pi_{n}$ is
approximately $\frac{1}{2}n^{2}$, and thus the cost of evaluating
$\mathcal{B}_{n}$ is only approximately 4 times the dimension of $\Pi_{n}$.
Qualitatively this is the same as in the univariate case. To evaluate a
polynomial
\begin{equation}
p\left(  x,y\right)  =\sum_{j=0}^{n}\sum_{k=0}^{j}b_{j,k}Q_{j}^{k}\left(
x,y\right)  \label{PolyOrtho}%
\end{equation}
for which $\left\{  b_{j,k}\right\}  $ are given, we use
\[
2\left(  n^{2}+3n-3\right)  +\left(  n+1\right)  \left(  n+2\right)
\approx3n^{2}%
\]
arithmetic operations, approximately 6 times the dimension $N_{n}$ of $\Pi
_{n}$.

There are other known choices of an orthonormal basis for $\Pi_{n}$; see Dunkl
and Xu \cite[\S 2.3.2]{DX} and Xu \cite[\S 1.2]{Xu2004}. In a number of
previous papers (see \cite{ach2008}, \cite{ah2010}, \cite{ahc2009},
\cite{ahc2013}) we have used the `ridge polynomials' of \cite{Loga}, in large
part because of their simple analytic form \ that is based on Chebyshev
polynomials of the second kind. However, we have calculated experimentally the
matrices $A_{i,n}$ and have found them to be dense for low order cases,
leading us to believe the same is true for larger values of $n$. For that
reason, solving the triple recursion relation (\ref{triple}) would be much
more costly than $\mathcal{O}\left(  n\right)  $ operations, making the choice
(\ref{OrthoBasis}) preferable in computational cost. As a particular example
of the lack of sparsity in the coefficient matrices $\left\{  A_{n,i}\right\}
$ for the ridge polynomials,
\[
A_{2,1}=\left[
\begin{array}
[c]{cccc}%
\frac{1}{2} & 0 & 0 & 0\smallskip\\
0 & \frac{\sqrt{2}}{8}+\frac{\sqrt{6}}{12} & -\frac{\sqrt{3}}{12} &
-\frac{\sqrt{2}}{8}+\frac{\sqrt{6}}{12}\smallskip\\
0 & \frac{\sqrt{2}}{8}-\frac{\sqrt{6}}{12} & \frac{\sqrt{3}}{12} &
-\frac{\sqrt{2}}{8}-\frac{\sqrt{6}}{12}%
\end{array}
\right]
\]%
\[
A_{2,2}=\left[
\begin{array}
[c]{cccc}%
0 & \frac{\sqrt{2}}{6} & -\frac{1}{6} & \frac{\sqrt{2}}{6}\smallskip\\
-\frac{\sqrt{3}}{12} & \frac{\sqrt{2}}{24}+\frac{\sqrt{6}}{12} & \frac{1}{3} &
\frac{\sqrt{2}}{24}-\frac{\sqrt{6}}{12}\smallskip\\
\frac{\sqrt{3}}{12} & \frac{\sqrt{2}}{24}-\frac{\sqrt{6}}{12} & \frac{1}{3} &
\frac{\sqrt{2}}{24}+\frac{\sqrt{6}}{12}%
\end{array}
\right]
\]

\subsection{Evaluating derivatives}

First derivatives of the orthonormal polynomials are required when
implementing the spectral methods of \cite{ach2008}, \cite{ah2010},
\cite{ahc2009}, \cite{ahc2013}). From (\ref{poly1}), (\ref{poly2}),%
\[%
\begin{array}
[c]{ll}%
\dfrac{\partial Q_{0}^{0}}{\partial x_{1}}=0, & \dfrac{\partial Q_{0}^{0}%
}{\partial x_{2}}=0,\medskip\\
\dfrac{\partial Q_{1}^{0}}{\partial x_{1}}=\dfrac{2}{\sqrt{\pi}}, &
\dfrac{\partial Q_{1}^{0}}{\partial x_{2}}=0\medskip\\
\dfrac{\partial Q_{1}^{1}}{\partial x_{1}}=0, & \dfrac{\partial Q_{1}^{1}%
}{\partial x_{2}}=\dfrac{2}{\sqrt{\pi}}%
\end{array}
\]
To obtain the first derivatives of the higher degree polynomials, we
differentiate the triple recursion relations of (\ref{q1})-(\ref{q3}). In
particular,%
\begin{equation}%
\begin{array}
[c]{l}%
\dfrac{\partial Q_{n+1}^{i}}{\partial x_{1}}=\dfrac{1}{a_{i,n}}\left\{
Q_{n}^{i}+x_{1}\dfrac{\partial Q_{n}^{i}}{\partial x_{1}}-a_{i,n-1}%
\dfrac{\partial Q_{n-1}^{i}}{\partial x_{1}}\right\}  ,\quad\quad
i=0,1,\dots,n-1\medskip\\
\dfrac{\partial Q_{n+1}^{n}}{\partial x_{1}}=\dfrac{1}{a_{i,n}}\left\{
Q_{n}^{n}+x_{1}\dfrac{\partial Q_{n}^{n}}{\partial x_{1}}\right\}  \medskip\\
\dfrac{\partial Q_{n+1}^{n+1}}{\partial x_{1}}=\dfrac{1}{d_{n,n}}\left\{
x_{2}\dfrac{\partial Q_{n}^{n}}{\partial x_{1}}-c_{n,n}\dfrac{\partial
Q_{n+1}^{n-1}}{\partial x_{1}}-d_{n-1,n-1}\dfrac{\partial Q_{n-1}^{n-1}%
}{\partial x_{1}}\right\}
\end{array}
\label{deriv1}%
\end{equation}%
\begin{equation}%
\begin{array}
[c]{l}%
\dfrac{\partial Q_{n+1}^{i}}{\partial x_{2}}=\dfrac{1}{a_{i,n}}\left\{
x_{1}\dfrac{\partial Q_{n}^{i}}{\partial x_{2}}-a_{i,n-1}\dfrac{\partial
Q_{n-1}^{i}}{\partial x_{2}}\right\}  ,\quad\quad i=0,1,\dots,n-1\medskip\\
\dfrac{\partial Q_{n+1}^{n}}{\partial x_{2}}=\dfrac{x_{1}}{a_{i,n}}%
\dfrac{\partial Q_{n}^{n}}{\partial x_{2}}\medskip\\
\dfrac{\partial Q_{n+1}^{n+1}}{\partial x_{2}}=\dfrac{1}{d_{n,n}}\left\{
Q_{n}^{n}+x_{2}\dfrac{\partial Q_{n}^{n}}{\partial x_{2}}-c_{n,n}%
\dfrac{\partial Q_{n+1}^{n-1}}{\partial x_{2}}-d_{n-1,n-1}\dfrac{\partial
Q_{n-1}^{n-1}}{\partial x_{2}}\right\}
\end{array}
\label{deriv2}%
\end{equation}

\section{Least squares approximation\label{sec3}}

When given a function $f\in C\left(  \mathbb{B}_{2}\right)  $, we are
interested in obtaining the least squares approximation to $f$ from the
polynomial subspace $\Pi_{n}$. When given the basis $\mathcal{B}_{n}$, this
approximation is given by the truncated Fourier expansion%
\begin{equation}
\mathcal{Q}_{n}f\left(  x,y\right)  \equiv P_{n}\left(  x,y\right)
=\sum_{j=0}^{n}\sum_{k=0}^{j}\left(  f,Q_{m}^{k}\right)  Q_{m}^{k}\left(
x,y\right)  . \label{LstSqApprox}%
\end{equation}
The linear operator $\mathcal{Q}_{n}$ is the orthogonal projection of
$L^{2}\left(  \mathbb{B}_{2}\right)  $ onto $\Pi_{n}$. As an operator on
$L^{2}\left(  \mathbb{B}_{2}\right)  $, it has norm 1. As an operator on
$C\left(  \mathbb{B}_{2}\right)  $ with the uniform norm $\left\Vert
\cdot\right\Vert _{\infty}$, $\mathcal{Q}_{n}$ has norm $\mathcal{O}\left(
n\right)  $; see \cite{XuB}.

The Fourier coefficients $\left(  f,Q_{m}^{k}\right)  $ must be evaluated
numerically, and we review a standard quadrature scheme to do so. Use the
formula%
\begin{equation}
\int_{\mathbb{B}_{2}}g(x,y)\,dx\,dy\approx\sum_{l=0}^{q}\sum_{m=0}%
^{2q}g\left(  r_{l},\frac{2\pi\,m}{2q+1}\right)  \omega_{l}\frac{2\pi}%
{2q+1}r_{l} \label{quad}%
\end{equation}
Here the numbers $r_{l}$ and $\omega_{l}$ are the nodes and weights,
respectively, of the $\left(  q+1\right)  $-point Gauss-Legendre quadrature
formula on $[0,1]$. Note that
\[
\int_{0}^{1}p(x)dx=\sum_{l=0}^{q}p(r_{l})\omega_{l},
\]
for all single-variable polynomials $p(x)$ with $\deg\left(  p\right)
\leq2q+1 $. The formula (\ref{quad}) uses the trapezoidal rule with $2q+1$
subdivisions for the integration over $\mathbb{B}_{2}$ in the azimuthal
variable. This quadrature is exact for all polynomials $g\in\Pi_{2q}$. For
functions $f,g\in C\left(  \mathbb{B}_{2}\right)  $, let $\left(  f,g\right)
_{q}$ denote the approximation of $\left(  f,g\right)  $ by the scheme
(\ref{quad}).

Our discrete approximation to (\ref{LstSqApprox}) is
\begin{equation}
\widetilde{P}_{n,q}\left(  x,y\right)  =\sum_{j=0}^{n}\sum_{k=0}^{j}\left(
f,Q_{m}^{k}\right)  _{q}Q_{m}^{k}\left(  x,y\right)  \label{DscLstSq}%
\end{equation}
When $q=n$, this approximation is known as the `discrete orthogonal projection
of $f$ onto $\Pi_{n}$', `hyperinterpolation of $f$ by $\Pi_{n}$', or the
`discrete least squares approximation'. We denote it by
\[
\widetilde{\mathcal{Q}}_{n}f\left(  x,y\right)  \equiv\widetilde{P}%
_{n,n}\left(  x,y\right)  \equiv\widetilde{P}_{n}\left(  x,y\right)
\]
In applying this numerical integration to the coefficients $\left(
f,Q_{m}^{k}\right)  $, we always require $q\geq n$ in order to force the
formula (\ref{LstSqApprox}) to reproduce all polynomials $f\in\Pi_{n}$. With
this requirement,
\[
f\in\Pi_{n}\quad\Rightarrow\quad\widetilde{P}_{n,q}=f.
\]
The operator $\widetilde{\mathcal{Q}}_{n}$ is a discrete orthogonal projection
of $C\left(  \mathbb{B}_{2}\right)  $ onto $\Pi_{n}$. For this specific case
of approximation over $\mathbb{B}_{2}$, see the discussion in \cite{hac}. In
particular,%
\[
\left\Vert \widetilde{\mathcal{Q}}_{n}\right\Vert _{C\rightarrow
C}=\mathcal{O}\left(  n\log n\right)  .
\]

\subsection{Cost of the discrete least squares approximation}

The main computational cost in (\ref{DscLstSq}) is the evaluation of the
coefficients $\left\{  \left(  f,Q_{m}^{k}\right)  _{q}\right\}  $. We begin
with the evaluation of the basis $\mathcal{B}_{n}$ at the \ points used in
(\ref{quad}), of which there are
\[
\left(  q+1\right)  \left(  2q+1\right)  .
\]
The cost to evaluate $\mathcal{B}_{n}$ will be
\begin{equation}
2\left(  n^{2}+3n-3\right)  \times\left(  q+1\right)  \left(  2q+1\right)
\approx4n^{2}q^{2} \label{cost1}%
\end{equation}
arithmetic operations.\ For comparison, recall that the dimension of $\Pi_{n}$
is approximately $\frac{1}{2}n^{2}$. The evaluation of the function $f$ at
these same nodes is
\begin{equation}
\left(  q+1\right)  \left(  2q+1\right)  N_{f}, \label{cost2}%
\end{equation}
with $N_{f}$ the cost of an individual evaluation of the function $f$. The
subsequent evaluations of the coefficients $\left\{  \left(  f,Q_{m}%
^{k}\right)  _{q}\right\}  $ involves an additional%
\begin{equation}
\frac{1}{2}\left(  n+1\right)  \left(  n+2\right)  \times\left(  q+1\right)
\left(  2q+1\right)  \label{cost3}%
\end{equation}
arithmetic operations. Having the coefficients $\left\{  \left(  f,Q_{m}%
^{k}\right)  _{q}\right\}  $, the polynomial (\ref{DscLstSq}) then requires
\begin{equation}
4\left(  n^{2}+3n-3\right)  \label{cost4}%
\end{equation}
arithmetic operations for each evaluation point $\left(  x,y\right)  $.

In the case $q=n$, the evaluation of $\widetilde{\mathcal{Q}}_{n}f$ is
dominated by (\ref{cost1}) and (\ref{cost3}), approximately $5n^{4}$
arithmetic operations. If we then evaluate $\widetilde{\mathcal{Q}}%
_{n}f\left(  x,y\right)  $ at the points used in the quadrature formula
(\ref{quad}), then the cost is an additional $8n^{4}$ operations, approximately.

\subsection{Convergence of least squares approximation}

Because the polynomials are dense in $L^{2}\left(  \mathbb{B}_{2}\right)  $,
we have
\[
\left\Vert f-P_{n}\right\Vert _{L^{2}}\rightarrow0\quad\text{as}\quad
n\rightarrow\infty.
\]
For convergence in $C\left(  \mathbb{B}_{2}\right)  $, we refer to the
presentation in \cite[\S 4.3.3, \S 5.7.1]{AtkinsonHan}. In particular,%
\begin{align}
\left\Vert f-\mathcal{Q}_{n}f\right\Vert _{\infty}  &  \leq\left(
1+\left\Vert \mathcal{Q}_{n}\right\Vert \right)  E_{n,\infty}\left(  f\right)
\medskip\label{Bound1}\\
\left\Vert f-\widetilde{\mathcal{Q}}_{n}f\right\Vert _{\infty}  &  \leq\left(
1+\left\Vert \widetilde{\mathcal{Q}}_{n}\right\Vert \right)  E_{n,\infty
}\left(  f\right)  \label{Bound2}%
\end{align}
where%
\[
E_{n,\infty}\left(  f\right)  =\min_{f\in\Pi_{n}}\left\Vert f-p\right\Vert
_{\infty},
\]
the minimax error in the approximation of $f$ by polynomials from $\Pi_{n}$.

Let $f\in C^{k,\alpha}\left(  \mathbb{B}_{2}\right)  $, functions that are $k-
$times continuously differentiable and whose $k^{\text{th}}$ derivatives are
H\"{o}lder continuous with exponent $\alpha\in(0,1]$ Then%
\begin{align}
E_{n,\infty}\left(  f\right)   &  \leq\frac{c_{k,\alpha}\left(  f\right)
}{n^{k+\alpha}},\quad\quad n\geq1. \label{Bound3}%
\end{align}
Combining these results with (\ref{Bound1})-(\ref{Bound2}) gives uniform
convergence of both $\mathcal{Q}_{n}f$ and $\widetilde{\mathcal{Q}}_{n}f$ to
$f$ for all $f\in C^{k,\alpha}\left(  \mathbb{B}_{2}\right)  $ with $k\geq1$.

\section{Triple recursion relation over the unit ball\label{sec4}}

In this section we repeat for the three dimensional case some of the results
from the two dimensional case of Sections \ref{sec2} and \ref{sec3}. The
orthonormal polynomials in this case are again taken from \cite[Proposition
2.3.2]{DX}. Here we first derive the coefficients of the three term recursion
relation in (\ref{triple2}).

\subsection{The recursion coefficients and the three term recurrence}

The orthonormal polynomials for the three dimensional unit ball are given by
\begin{align}
Q_{n}^{j,k}(x,y,z)  &  =\frac{1}{h_{j,k}}C_{n-j-k}^{j+k+3/2}(x)(1-x^{2}%
)^{j/2}\nonumber\\
&  \cdot C_{j}^{k+1}(\frac{y}{\sqrt{1-x^{2}}})(1-x^{2}-y^{2})^{k/2}C_{k}%
^{1/2}(\frac{z}{\sqrt{1-x^{2}-y^{2}}}) \label{eq4001}%
\end{align}
where $j+k\leq n$, and $n\in\mathbb{N}$ is the degree of the polynomial
$Q_{n}^{j,k}$. The normalization constant $h_{j,k}$ will be derived further
below. We introduce the vector of all orthonormal polynomials $\mathbb{P}_{n}
$ of degree $n$:
\begin{equation}
\mathbb{P}_{n}=[Q_{n}^{0,0},\ldots,Q_{n}^{0,n},Q_{n}^{1,0},\ldots
,Q_{n}^{1,n-1},Q_{n}^{2,0},\dots,Q_{n}^{2,n-2},\dots,Q_{n}^{n,0}]^{T},\quad
n\geq0. \label{eq4002}%
\end{equation}
Here we have ${\binom{n+2}{2}}$ polynomials of degree $n$ and the space
$\Pi_{n}$ has dimension ${\binom{n+3}{3}}$, see \cite{DX}. In formula
(\ref{triple2}) we have matrices $A_{n,i}$, $i=1,2,3$, of dimension
${\binom{n+2}{2}}\times{\binom{n+3}{2}}$. First we derive the normalization
constant $h_{j,k}$ with a calculation which is typical for calculations
involved in the calculation of the coefficients of the matrices $A_{n,i}$. By
definition we have
\begin{align}
h_{j,k}^{2}  &  =\int_{-1}^{1}(C_{n-j-k}^{j+k+3/2}(x))^{2}(1-x^{2})^{j}%
\int_{-\sqrt{1-x^{2}}}^{\sqrt{1-x^{2}}}(C_{j}^{k+1}(\frac{y}{\sqrt{1-x^{2}}%
}))^{2}(1-x^{2}-y^{2})^{k}\nonumber\\
&  \cdot\int_{-\sqrt{1-x^{2}-y^{2}}}^{\sqrt{1-x^{2}-y^{2}}}(C_{k}^{1/2}%
(\frac{z}{\sqrt{1-x^{2}-y^{2}}}))^{2}dz\;dy\;dx \label{eq4004}%
\end{align}
Using the substitution
\begin{align*}
u:=  &  \frac{z}{\sqrt{1-x^{2}-y^{2}}}\\
dz  &  =\sqrt{1-x^{2}-y^{2}}\;du
\end{align*}
we get
\begin{align*}
h_{j,k}^{2}  &  =\int_{-1}^{1}(C_{n-j-k}^{j+k+3/2}(x))^{2}(1-x^{2})^{j}%
\int_{-\sqrt{1-x^{2}}}^{\sqrt{1-x^{2}}}(C_{j}^{k+1}(\frac{y}{\sqrt{1-x^{2}}%
}))^{2}(1-x^{2}-y^{2})^{k+1/2}\\
&  \cdot\int_{-1}^{1}(C_{k}^{1/2}(u))^{2}du\;dy\;dx\\
&  =N_{k}^{[1/2]}\int_{-1}^{1}(C_{n-j-k}^{j+k+3/2}(x))^{2}(1-x^{2})^{j}\\
&  \cdot\int_{-\sqrt{1-x^{2}}}^{\sqrt{1-x^{2}}}(C_{j}^{k+1}(\frac{y}%
{\sqrt{1-x^{2}}}))^{2}(1-x^{2}-y^{2})^{k+1/2}\;dy\;dx
\end{align*}
where we defined
\begin{align}
(N_{k}^{[\mu]})^{2}:=  &  \int_{-1}^{1}(C_{k}^{\mu}(x))^{2}(1-x^{2})^{\mu
-1/2}dx\nonumber\\
&  =\frac{\pi\Gamma(2\mu+k)}{2^{2\mu-1}k!(\mu+k)\Gamma^{2}(\mu)}
\label{eq4003}%
\end{align}
see \cite{AS}. Now we use the substitution
\begin{align*}
u:=  &  \frac{y}{\sqrt{1-x^{2}}}\\
dz  &  =\sqrt{1-x^{2}}\;du\\
(1-x^{2}-y^{2})  &  =(1-x^{2}-(1-x^{2})u^{2})\\
&  =(1-x^{2})(1-u^{2})
\end{align*}
to obtain
\begin{align*}
h_{j,k}^{2}  &  =N_{k}^{[1/2]}\int_{-1}^{1}(C_{n-j-k}^{j+k+3/2}(x))^{2}%
(1-x^{2})^{j+k+1}\\
&  \cdot\int_{-1}^{1}(C_{j}^{k+1}(u))^{2}(1-u^{2})^{k+1/2}du\;dx\\
&  =N_{k}^{[1/2]}N_{j}^{[k+1]}N_{n-j-k}^{[j+k+3/2]}%
\end{align*}
If we denote the coefficients of the matrices $A_{n,i}$ by $a_{j,k;j^{\prime
},k^{\prime}}^{[n,i]}$ $j+k\leq n$ and $j^{\prime}+k^{\prime}\leq n+1$ we get
\begin{align*}
a_{j,k;j^{\prime},k^{\prime}}^{[n,1]}  &  =\int_{\mathbb{B}_{3}}xQ_{n}%
^{j,k}(x,y,z)Q_{n+1}^{j^{\prime},k^{\prime}}(x,y,z)\;d(x,y,z)\\
a_{j,k;j^{\prime},k^{\prime}}^{[n,2]}  &  =\int_{\mathbb{B}_{3}}yQ_{n}%
^{j,k}(x,y,z)Q_{n+1}^{j^{\prime},k^{\prime}}(x,y,z)\;d(x,y,z)\\
a_{j,k;j^{\prime},k^{\prime}}^{[n,3]}  &  =\int_{\mathbb{B}_{3}}zQ_{n}%
^{j,k}(x,y,z)Q_{n+1}^{j^{\prime},k^{\prime}}(x,y,z)\;d(x,y,z)\\
&
\end{align*}
Each of the integrals can be written in the same way as the integral in
(\ref{eq4004}) and then the two above substitutions together with the
orthonormal property of the Gegenbauer polynomials allows us to calculate the
coefficients of $A_{n,i}$, $i=1,2,3$. Again we obtain very sparsely populated
matrices. Equation (\ref{triple2}) takes on the following form:
\begin{equation}
xQ_{n}^{j,k}=a_{j,k;j,k}^{[n,1]}Q_{n+1}^{j,k}+a_{j,k;j,k}^{[n-1,1]}%
Q_{n-1}^{j,k},:j+k\leq n \label{eq4005}%
\end{equation}
where
\begin{equation}
a_{j,k;j,k}^{[n,1]}=\frac{1}{2}\Big(\frac{(j+k+n+3)(n+1-j-k)}{(n+5/2)(n+3/2)}%
\Big)^{1/2} \label{eq4006}%
\end{equation}
and the term $a_{j,k;j,k}^{[n-1,1]}$ has to be replaced by $0$ if $j+k=n$.
Furthermore we get
\begin{align}
yQ_{n}^{j,k}  &  =a_{j,k;j+1,k}^{[n,2]}Q_{n+1}^{j+1,k}+a_{j,k;j-1,k}%
^{[n,2]}Q_{n+1}^{j-1,k}\nonumber\\
&  +a_{j+1,k;j,k}^{[n+1,2]}Q_{n-1}^{j+1,k}+a_{j-1,k;j,k}^{[n-1,2]}%
Q_{n-1}^{j-1,k},:j+k\leq n \label{eq4007}%
\end{align}
where the terms of the matrix $A_{n-1,2}$ and $A_{n,2}$ have to substituted by
zero if $j-1+k<0$ or $j+1+k>n-1$ in the case of $A_{n-1,2}$. Here
\begin{align}
a_{j,k;j+1,k}^{[n,2]}  &  =\frac{1}{4}\Big(\frac
{(j+2k+2)(j+1)(j+k+n+4)(j+k+n+3)}{(j+k+1)(j+k+2)(n+5/2)(n+3/2)}\Big)^{1/2}%
\label{eq4008}\\
a_{j,k;j-1,k}^{[n,2]}  &  =-\frac{1}{4}\Big(\frac{j(j+2k+1)(n+2-j-k)(n+1-j-k)}%
{(j+k+1)(j+k)(n+3/2)(n+5/2)}\Big)^{1/2} \label{eq4009}%
\end{align}
Finally we get
\begin{align}
zQ_{n}^{j,k}  &  =a_{j,k;j,k-1}^{[n,3]}Q_{n+1}^{j,k-1}+a_{j,k;j+2,k-1}%
^{[n,3]}Q_{n+1}^{j+2,k-1}+a_{j,k;j,k+1}^{[n,3]}Q_{n+1}^{j,k+1}\nonumber\\
&  +a_{j,k;j-2,k+1}^{[n,3]}Q_{n+1}^{j-2,k+1}+a_{j+2,k-1;j,k}^{[n-1,3]}%
Q_{n-1}^{j+2,k-1}+a_{j,k-1;j,k}^{[n-1,3]}Q_{n-1}^{j,k-1}\nonumber\\
&  +a_{j-2,k+1;j,k}^{[n-1,3]}Q_{n-1}^{j-2,k+1}+a_{j,k+1;j,k}^{[n-1,3]}%
Q_{n-1}^{j,k+1},:j+k\leq n \label{eq4010}%
\end{align}
where again the terms have to be replaced by zero if the indices are out of
the range of the corresponding matrix. Here
\begin{align}
&  a_{j,k;j,k-1}^{[n,3]}\nonumber\\
&  =-\frac{k}{8}\Big(\frac{(j+2k+1)(j+2k)(n+2-j-k)(n+1-j-k)}%
{(k+1/2)(k-1/2)(j+k+1)(j+k)(n+3/2)(n+5/2)}\Big)^{1/2}\label{eq4011}\\
&  a_{j,k;j+2,k-1}^{[n,3]}\nonumber\\
&  =-\frac{k}{8}\Big(\frac{(j+2)(j+1)(j+k+n+4)(j+k+n+3)}%
{(k+1/2)(k-1/2)(j+k+1)(j+k+2)(n+3/2)(n+5/2)}\Big)^{1/2}\label{eq4012}\\
&  a_{j,k;j,k+1}^{[n,3]}\nonumber\\
&  =\frac{k+1}{8}\Big(\frac{(j+2k+3)(j+2k+2)(j+k+n+4)(j+k+n+3)}%
{(k+1/2)(k+3/2)(j+k+1)(j+k+2)(n+3/2)(n+5/2)}\Big)^{1/2}\label{eq4013}\\
&  a_{j,k;j-2,k+1}^{[n,3]}\nonumber\\
&  =\frac{k+1}{8}\Big(\frac{(n+2-j-k)(n+1-j-k)j(j-1)}%
{(k+1/2)(k+3/2)(j+k)(j+k+1)(n+3/2)(n+5/2)}\Big)^{1/2} \label{eq4014}%
\end{align}
The equations (\ref{eq4005}), (\ref{eq4007}), and (\ref{eq4010}) allow the
calculation of all $Q_{n+1}^{j,k}$ in the following way. For $j+k\leq n$ we
can use (\ref{eq4005}) and solve for $Q_{n+1}^{j,k}$:
\begin{equation}
Q_{n+1}^{j,k}=\frac{xQ_{n}^{j,k}-a_{j,k;j,k}^{[n-1,1]}Q_{n-1}^{j,k}%
}{a_{j,k;j,k}^{[n,1]}} \label{eq4015}%
\end{equation}
Then we use (\ref{eq4007}) for the calculation of $Q_{n+1}^{j+1,n-j}$,
$j=0,\ldots,n$:
\begin{align}
Q_{n+1}^{j+1,n-j}  &  =\Big(yQ_{n}^{j,n-j}-a_{j,n-j;j-1,n-j}^{[n,2]}%
Q_{n+1}^{j-1,n-j}\nonumber\\
&  -a_{j-1,n-j;j,n-j}^{[n-1,2]}Q_{n-1}^{j-1,n-j}\Big)\Big/a_{j,n-j;j+1,n-j}%
^{[n,2]} \label{eq4016}%
\end{align}
Finally (\ref{eq4010}) allows us to calculate $Q_{n+1}^{0,n+1}$%

\begin{align}
Q^{0,n+1}_{n+1}  &  = \Big( zQ^{0,n}_{n} - a^{[n,3]}_{0,n;0,n-1}%
Q^{0,n-1}_{n+1}- a^{[n,3]}_{0,n;2,n-1}Q^{2,n-1}_{n+1}\nonumber\\
&  -a^{[n-1,3]}_{0,n-1;0,n}Q^{0,n-1}_{n-1} \Big)\Big/ a^{[n,3]}_{0,n;0,n+1}
\label{eq4017}%
\end{align}
By taking partial derivatives in equation (\ref{eq4015})--(\ref{eq4017}) we
are able to derive recursion formulas for the partial derivatives of the
orthonormal polynomials as in (\ref{deriv1})--(\ref{deriv2}).

\subsection{Least square approximation}

Similar to Section 3, the least square approximation in $L^{2}(\mathbb{B}%
_{3})$ for a function $f\in L^{2}(\mathbb{B}_{3})$ is given by
\begin{equation}
\mathcal{Q}_{n}f(x,y,z)=P_{n}(x,y,z)=\sum_{m=0}^{n}\sum_{j+k\leq m}%
(f,Q_{m}^{j,k})Q_{m}^{j,k}(x,y,z) \label{eq4019}%
\end{equation}
where the inner product is given by
\begin{equation}
(f,Q_{m}^{j,k})=\int_{\mathbb{B}_{3}}f(x,y,z)Q_{m}^{j,k}(x,y,z)\;d(x,y,z)
\label{eq4020}%
\end{equation}
For practical calculations we have to replace the integral in (\ref{eq4020})
by a quadrature rule for $f\in C(\mathbb{B}_{3})$. One choice is to use a
quadrature rule which will integrate polynomials of degree smaller or equal to
$2n$ exactly, so we have
\begin{equation}
\mathcal{Q}_{n}p(x,y,z)=p(x,y,z),\quad\forall p\in\Pi_{n} \label{eq4021}%
\end{equation}
We will use
\begin{align}
\int_{\mathbb{B}_{3}}g(x,y,z)\,d(x,y,z)  &  =\int_{0}^{1}\int_{0}^{2\pi}%
\int_{0}^{\pi}\widetilde{g}(r,\theta,\phi)\,r^{2}\sin(\phi)\,d\phi
\,d\theta\,dr\approx Q_{q}[g]\medskip\nonumber\\
Q_{q}[g]  &  :=\sum_{i=1}^{2q}\sum_{j=1}^{q}\sum_{k=1}^{q}\frac{\pi}%
{q}\,\omega_{j}\,\nu_{k}\widetilde{g}\left(  \frac{\zeta_{k}+1}{2},\frac
{\pi\;i}{2q},\arccos(\xi_{j})\right)  \label{eq4022}%
\end{align}
$q>n$. Here $\widetilde{g}(r,\theta,\phi)=g(x,y,z)$ is the representation of
$g$ in spherical coordinates. For the $\theta$ integration we use the
trapezoidal rule, because the function is $2\pi-$periodic in $\theta$. For the
$r$ direction we use the transformation
\begin{align*}
\int_{0}^{1}r^{2}v(r)\;dr  &  =\int_{-1}^{1}\left(  \frac{t+1}{2}\right)
^{2}v\left(  \frac{t+1}{2}\right)  \frac{dt}{2}\medskip\\
&  =\frac{1}{8}\int_{-1}^{1}(t+1)^{2}v\left(  \frac{t+1}{2}\right)
\;dt\medskip\\
&  \approx\sum_{k=1}^{q}\underset{_{=:\nu_{k}}}{\underbrace{\frac{1}{8}\nu
_{k}^{\prime}}}v\left(  \frac{\zeta_{k}+1}{2}\right)
\end{align*}
where the $\nu_{k}^{\prime}$ and $\zeta_{k}$ are the weights and the nodes of
the Gauss quadrature with $q$ nodes on $[-1,1]$ with respect to the inner
product
\[
(v,w)=\int_{-1}^{1}(1+t)^{2}v(t)w(t)\,dt
\]
The weights and nodes also depend on $q$ but we omit this index. For the
$\phi$ direction we use the transformation
\begin{align*}
\int_{0}^{\pi}\sin(\phi)v(\phi)\,d\phi &  =\int_{-1}^{1}v(\arccos
(\phi))\,d\phi\medskip\\
&  \approx\sum_{j=1}^{q}\omega_{j}v(\arccos(\xi_{j}))
\end{align*}
where the $\omega_{j}$ and $\xi_{j}$ are the nodes and weights for the
Gauss--Legendre quadrature on $[-1,1]$. This quadrature rule has been used in
our earlier articles, see \cite{ach2008}. For more information on this
quadrature rule on the unit ball in $\mathbb{R}^{3}$, see \cite{stroud}. For
the complexity estimation in the next section we will assume that we use the
smallest possible $q$ to satisfy (\ref{eq4021}) which is $q=n+1$. Although a
little bit larger values of $q$ might improve the approximation property of
(\ref{eq4019}) in practice. With this value of $q$ the quadrature formula
(\ref{eq4022}) uses $2\left(  n+1\right)  ^{3}=2n^{3}+\mathcal{O}\left(
n^{2}\right)  $ points in the unit ball $\mathbb{B}_{3}$.

The discrete $L^{2}$ projection is now given by
\begin{equation}
\widetilde{\mathcal{Q}}_{n}f(x,y,z)=\widetilde{P}_{n}(x,y,z)=\sum_{m=0}%
^{n}\sum_{j+k\leq m}Q_{n}[f\cdot Q_{m}^{j,k}]Q_{m}^{j,k}(x,y,z) \label{eq4023}%
\end{equation}

Regarding the convergence of the convergence of $\mathcal{Q}_{n}f$ towards $f
$ in $L^{2}(\mathbb{B}_{3})$ and $L^{\infty}(\mathbb{B}_{3})$ we have similar
results to Section 3.2. Because the polynomials are dense we have convergence
in $L^{2}(\mathbb{B}_{3})$ and formulas (\ref{Bound1}) and (\ref{Bound2}) hold
as before, and the same is true for the estimate for $E_{n,\infty}(f)$ in
(\ref{Bound3}). But the Lebesgue constant for the projection $\mathcal{Q}_{n}$
in $L^{\infty}(\mathbb{B}_{3})$ is larger,
\begin{equation}
\Vert\mathcal{Q}_{n}\Vert_{C\mapsto C}=\mathcal{O}_{n\rightarrow\infty
}(n^{3/2}) \label{eq4024}%
\end{equation}
see \cite{XuB}. Together with (\ref{Bound3}) we obtain the convergence in
$C(\mathbb{B}_{3})$ for functions which are in $C^{1,\alpha}(\mathbb{B}_{3})$,
$\alpha>1/2$.

For the bound of $\| \widetilde{\mathcal{Q}}_{n} \|_{C \mapsto C}$ we can use
the same arguments as in (2.10)--(2.18) of our previous article \cite{hac}
together with the results about the reproducing kernel in \cite{XuB}. This
shows
\begin{align*}
\| \widetilde{\mathcal{Q}}_{n} \|_{C \mapsto C}  &  = \mathcal{O}%
_{n\rightarrow\infty}(n^{2})
\end{align*}
and proves the convergence of the discrete $L^{2}$ approximation in the
inifinity norm for functions which are in $C^{2,\alpha}(\mathbb{B}_{3})$,
$\alpha>0$.

\subsection{Computational cost}

First we give here a brief analysis of the computational cost to evaluate all
polynomials $Q_{m}^{j,k}$ in $\Pi_{n}$ at a given point. We assume again, that
all coefficients in (\ref{eq4015})--(\ref{eq4017}) have been calculated. If we
further assume that $Q_{m}^{j,k}$ and $Q_{m-1}^{j,k}$ have been calculated
then (\ref{eq4015}), for $j+k\leq m$, constitutes the dominant work for the
calculation of $Q_{m+1}^{j,k}$, $j+k\leq m+1$. To evaluate (\ref{eq4015}) for
$j+k\leq m$ requires $4{\binom{m+2}{2}}=2m^{2}+\mathcal{O}_{m\rightarrow
\infty}(m)$ arithmetic operations. The evaluation of (\ref{eq4016}) and
(\ref{eq4017}) will not change this asymptotic behavior. Adding these up for
$m=0,1,\ldots n-1$ leads to a total number of arithmetic operations given by
$\frac{2}{3}n^{3}+\mathcal{O}_{n\rightarrow\infty}(n^{2})$. If we further
consider the problem to evaluate the polynomial
\begin{equation}
p(x,y,z)=\sum_{m=0}^{n}\sum_{j+k\leq m}b_{m}^{j,k}Q_{m}^{j,k}(x,y,z)
\label{eq4018}%
\end{equation}
we have to add another $2\sum_{m=0}^{n}{\binom{m+2}{2}}=\frac{1}{3}%
n^{3}+\mathcal{O}_{n\rightarrow\infty}(n^{2})$ operations, which means that
the evaluation of (\ref{eq4018}) requires a total $n^{3}+\mathcal{O}%
_{n\rightarrow\infty}(n^{2})$ operations, if the recursion coefficients are
known. The set $\Pi_{n}$ has $\frac{n^{3}}{6}+\mathcal{O}_{n\rightarrow\infty
}(n^{2})$ elements, so about 6 operations are needed in average per basis
functions, exactly the same as in Section 2.

To calculate the discrete $L^{2}$ projection (\ref{eq4023}) we first need to
evaluate $f$ at the $\sim2n^{3}$ quadrature points of $Q_{n}$, this requires
an effort of $\sim2n^{3}N_{f}$, where $N_{f}$ again measures the cost of an
individual evaluation of $f$. Then we have to calculate all basis functions
$Q^{j,k}_{m}$ in $\Pi_{n}$ for all $2n^{3}$ points. This requires $\frac{4}%
{3}n^{6}+\mathcal{O}_{n\rightarrow\infty}(n^{5})$ operations. The calculation
of a single $Q_{n}[f\cdot Q^{j,k}_{m}]$ requires $6n^{3}$ operations and we
have to do this for all ${\binom{n+3 }{3}}$ basis functions of $\Pi_{n}$ which
results in an additional $n^{6}+\mathcal{O}_{n\rightarrow\infty}(n^{5})$
operations. If we assume that the $N_{f}$ is less than $\mathcal{O}(n^{3})$ we
see that the evaluation of the discrete inner products $Q_{n}[f\cdot
Q^{j,k}_{m}]$ is the dominant term and the complexity of the calculation of
(\ref{eq4023}) is given by $\frac{7}{3}n^{6}+\mathcal{O}_{n\rightarrow\infty
}(n^{5})$.

\section{Numerical examples and \textsc{Matlab} programs\label{sec5}}

We present \textsc{Matlab} programs for using orthonormal polynomials over the
unit disk. We compute the coefficients $\left\{  a_{i,n},c_{i,n}%
,d_{i,n}\right\}  $, the basis $\mathcal{B}_{n}$, and the discrete least
squares approximation (\ref{DscLstSq}) with $q\geq n$. The program
\texttt{TripleRecurCoeff} is used to produce the needed coefficients $\left\{
a_{i,n},c_{i,n},d_{i,n}\right\}  $, the program \texttt{EvalOrthoPolys} is
used to evaluate the polynomials in the basis $\mathcal{B}_{n}$, and the
program \texttt{LeastSqCoeff}\ evaluates the coefficients in (\ref{DscLstSq}).
The program \texttt{EvalLstSq} is used to evaluate $\widetilde{P}_{n,q}\left(
x,y\right)  $ at a selected set of nodes in $\mathbb{B}_{2}$; it also
evaluates the error and produces various graphs of the error as the degree $n$
is increased. The program \texttt{Test\_EvalLstSq} is used to test the
programs just listed.
\begin{figure}[tb]%
\centering
\includegraphics[
height=3in,
width=3.9998in
]%
{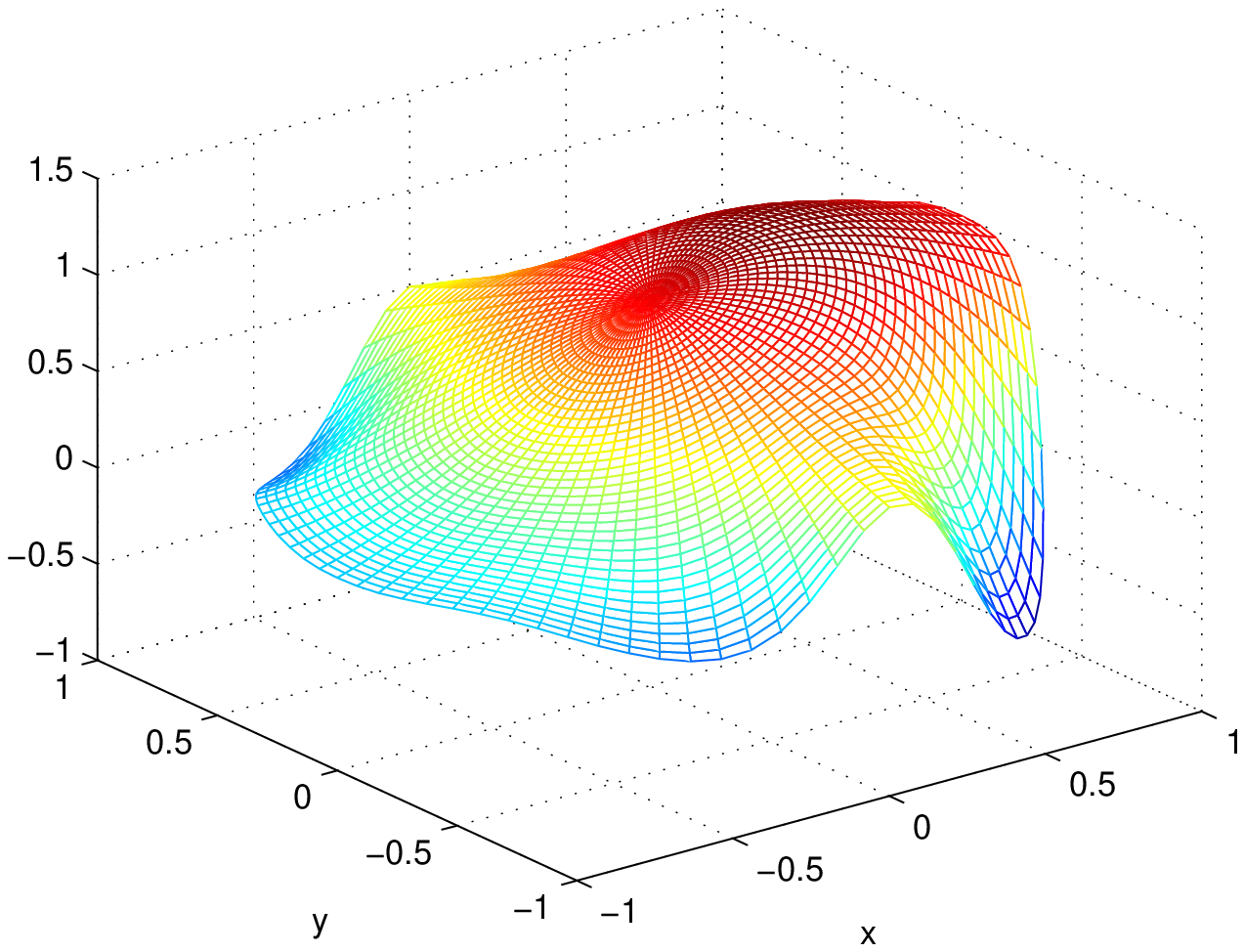}%
\caption{The approximation $\protect\widetilde{P}_{n,q}\left(  x,y\right)  $
for (\ref{TestFcn}), with $n=30$ and $q=40$}%
\label{approxdeg30}%
\end{figure}

\begin{figure}[tb]%
\centering
\includegraphics[
height=3in,
width=3.9998in
]%
{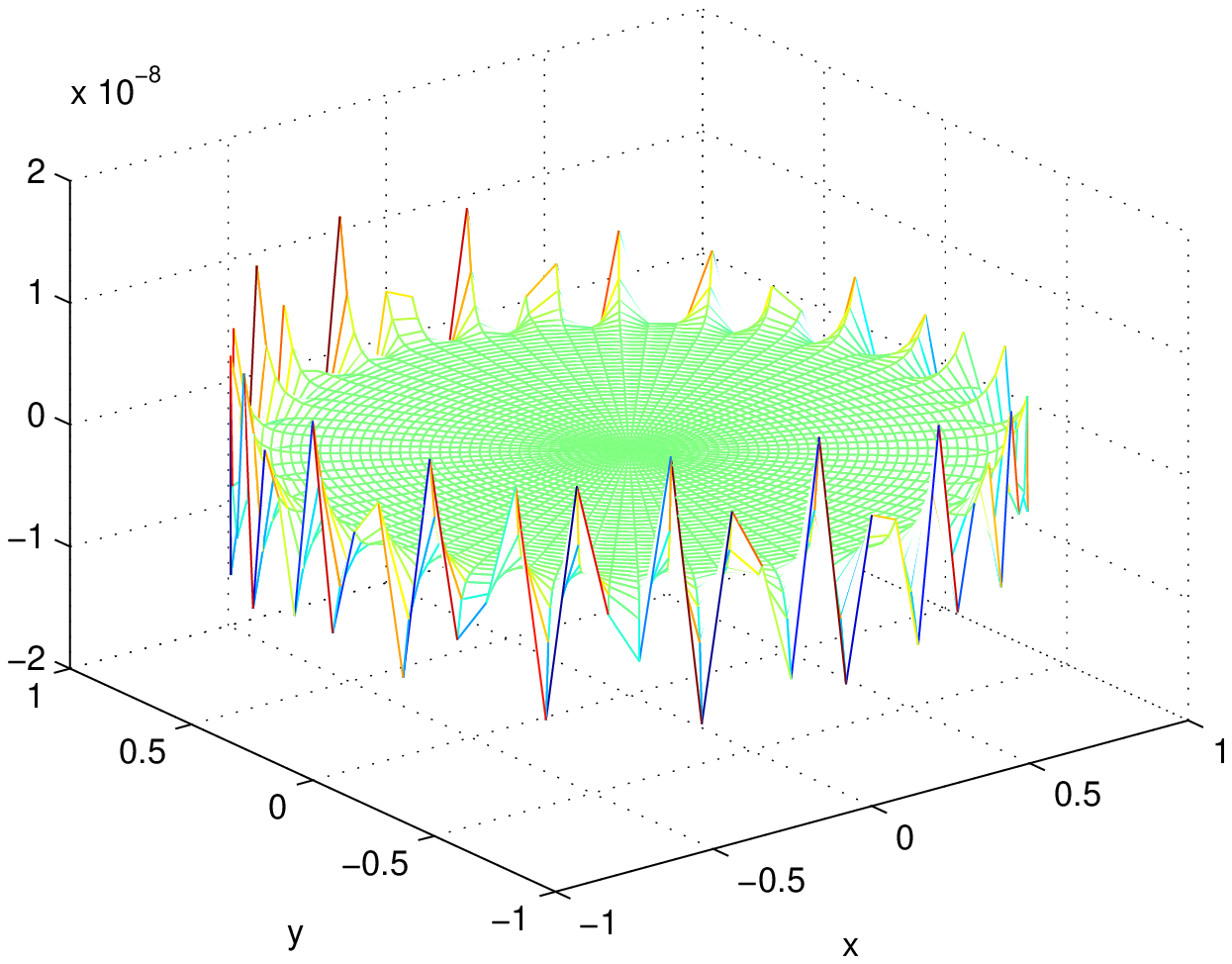}%
\caption{The error $f-$ $\protect\widetilde{P}_{n,q}$ for (\ref{TestFcn}),
with $n=30$ and $q=40$}%
\label{errordeg30}%
\end{figure}

\begin{figure}[tb]%
\centering
\includegraphics[
height=3in,
width=3.9998in
]%
{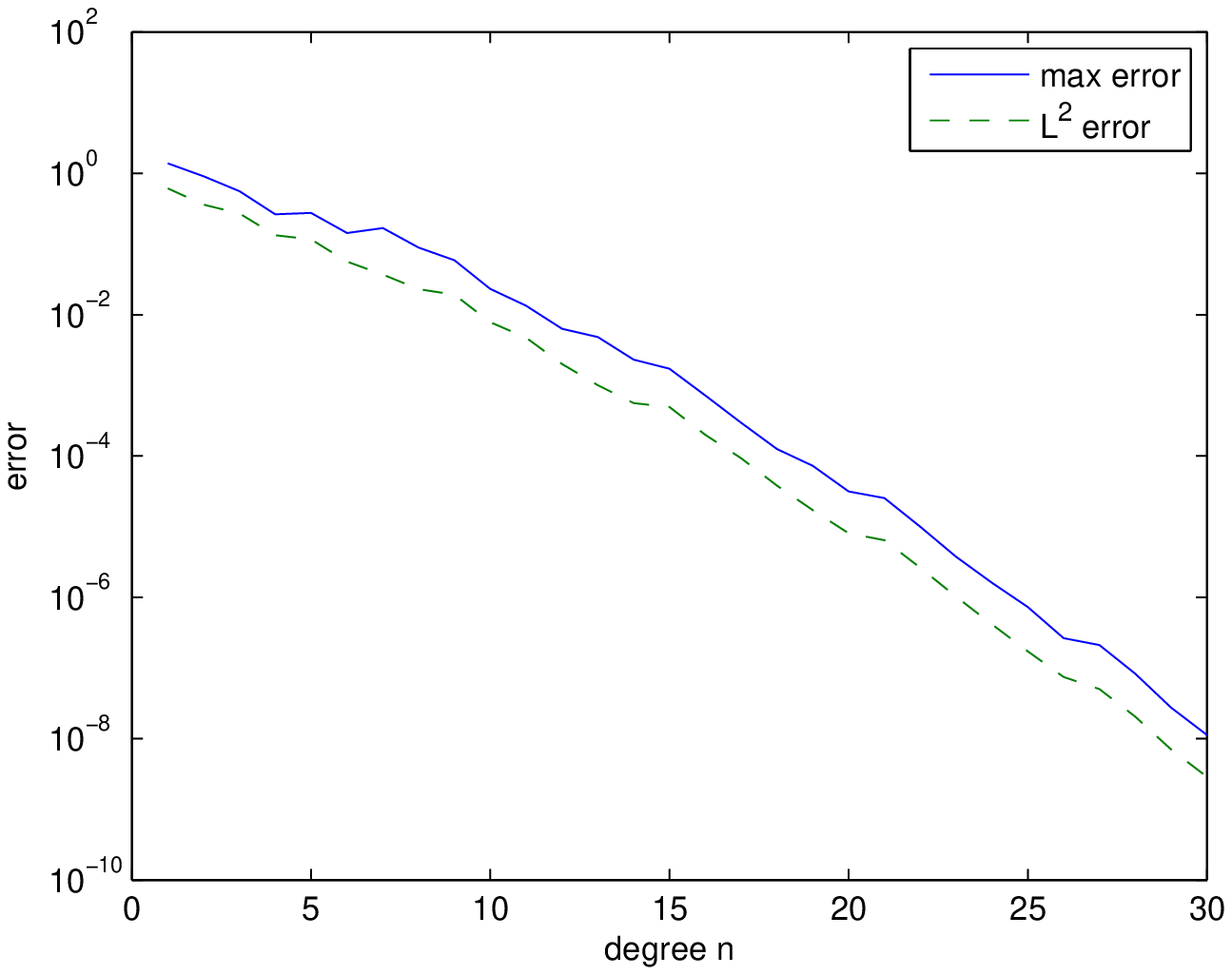}%
\caption{The error $f-$ $\protect\widetilde{P}_{n,q}$ for (\ref{TestFcn}) with
$q=40$}%
\label{degvserror}%
\end{figure}

Consider the function%
\begin{equation}
f\left(  x,y\right)  =\frac{1+x}{1+x^{2}+y^{2}}\cos\left(  6xy^{2}\right)
\label{TestFcn}%
\end{equation}
This was approximated using \texttt{Test\_EvalLstSq} for degrees 1 through 30.
Figure \ref{approxdeg30} shows $\widetilde{P}_{30,40}$ and Figure
\ref{errordeg30} shows its error. The error as it varies with the degree $n$
is shown in Figure \ref{degvserror}. This last graph suggests an exponential
rate of convergence for $\widetilde{P}_{n,q}$ to $f$.%

We have found often that the error $f\left(  x,y\right)  -\widetilde{P}%
_{n,q}\left(  x,y\right)  $ is slightly smaller than that of $f\left(
x,y\right)  -\widetilde{P}_{n,n}\left(  x,y\right)  $ if $q$ is taken a small
amount larger than $n$, say $q=n+5$. However, the qualitative behaviour shown
in Figure \ref{degvserror} is still valid for $f-\widetilde{P}_{n,n}$.

\subsection{Additional comments}

These programs can also be used for constructing approximations over other
planar regions $\Omega$. For example, the mapping
\[
\left(  x,y\right)  \mapsto\left(  \xi,\eta\right)  =\left(  ax,by\right)
,\quad\quad\left(  x,y\right)  \in\mathbb{B}_{2},
\]
with $a,b>0$, can be used to create polynomial approximations to a function
defined over the ellipse%
\[
\left(  \frac{\xi}{a}\right)  ^{2}+\left(  \frac{\eta}{b}\right)  ^{2}\leq1.
\]
If polynomials are not required, only an approximating function, then
mappings
\[
\left(  x,y\right)  \mapsto\left(  \xi,\eta\right)  =\Phi\left(  x,y\right)
,\quad\quad\left(  x,y\right)  \in\mathbb{B}_{2}%
\]
with $\Phi$ a 1-1 mapping can be used to convert an approximation problem over
a planar region $\Omega$ to one over $\mathbb{B}_{2}$. The construction of
such mappings $\Phi$ is explored in \cite{ah2011}.


\begin{thebibliography}{99}                                                                                               %


\bibitem {AS}M. Abramowitz, I. Stegun. \emph{Handbook of Mathematical
Functions}, Dover, 1965.

\bibitem {ach2008}K. Atkinson, D. Chien, and O. Hansen. A spectral method for
elliptic equations: The Dirichlet problem, \emph{Advances in Computational
Mathematics}\textit{, }\textbf{33} (2010), pp. 169-189.

\bibitem {AtkinsonHan}K. Atkinson and Weimin Han. \emph{Spherical Harmonics
and Approximations on the Unit Sphere : An Introduction}, Lecture Notes in
Mathematics \#2044, Springer-Verlag, New York, 2012.

\bibitem {ah2010}K. Atkinson and O. Hansen. A spectral method for the
eigenvalue problem for elliptic equations, \emph{Electronic Transactions on
Numerical Analysis}\textit{\ }\textbf{37} (2010), pp. 386-412.

\bibitem {ah2011}K. Atkinson and O. Hansen. Creating domain mappings,
\emph{Electronic Transactions on Numerical Analysis} \textbf{39} (2012), pp. 202-230.

\bibitem {ahc2009}K. Atkinson, O. Hansen, and D. Chien. A spectral method for
elliptic equations: The Neumann problem, \emph{Advances in Computational
Mathematics}\textit{\ }\textbf{34} (2011), pp. 295-317.

\bibitem {ahc2013}K. Atkinson,O. Hansen, and D. Chien. \ A spectral method for
parabolic differential equations, \emph{Numerical Algorithms}, DOI
10.1007/s11075-012-9620-8, to appear. A preliminary version is available at http://arxiv.org/abs/1203.6709

\bibitem {DX}C. Dunkl and Y. Xu. \emph{Orthogonal Polynomials of Several
Variables}, Cambridge Univ. Press, Cambridge, 2001.

\bibitem {hac}O. Hansen, K. Atkinson, and D. Chien. On the norm of the
hyperinterpolation operator on the unit disk and its use for the solution of
the nonlinear Poisson equation, \emph{IMA\ J. Numerical Analysis} \textbf{29}
(2009, 257-283, DOI: 10.1093/imanum/drm052.

\bibitem {Loga}B. Logan. and L. Shepp. Optimal reconstruction of a function
from its projections, \emph{Duke Mathematical Journal} \textbf{42}, (1975), 645--659.

\bibitem {nist}F. Olver, D. Lozier, R. Boisvert, C. Clark. \emph{NIST Handbook
of Mathematical Functions}, Cambridge University Press, 2010.

\bibitem {stroud}A. Stroud. \emph{Approximate Calculation of Multiple
Integrals,} Prentice-Hall, Inc., Englewood Cliffs, N.J., 1971.

\bibitem {XuB}Xu, Y. Representation of reproducing kernels and the Lebesgue
constants on the ball, \textit{J. Approx. Theor.} \textbf{112} (2001), 295--310.

\bibitem {Xu2004}Yuan Xu, Lecture notes on orthogonal polynomials of several
variables, in \emph{Advances in the Theory of Special Functions \& Orthogonal
Polynomials}, Nova Sci. Pub., 2004, 135-188.
\end{thebibliography}
\end{document}